\documentclass[12pt]{article}
\usepackage{amssymb,amsmath,amsfonts,amsthm,bm,latexsym, epsfig,cite, psfrag,eepic,colordvi,amscd,graphics}
\usepackage{lineno}
\usepackage{mathrsfs}

\renewcommand{\vec}[1]{\mbox{\boldmath$#1$}}
\newtheorem{thm}{Theorem}[section]

\newtheorem{cor}[thm]{Corollary}
\newtheorem{defn}[thm]{Definition}
\newtheorem{lem}[thm]{Lemma}
\newtheorem{fact}[thm]{Fact}

\title{A refinement of the formula for $k$-ary trees and the Gould-Vandermonde's convolution}

\author{Ricky X. F. Chen\\
\small Center for Combinatorics, LPMC-TJKLC\\[-0.8ex]
\small Nankai University, Tianjin 300071, P. R. China\\[-0.8ex]
\small \texttt{ricky\_chen@mail.nankai.edu.cn}}
\date{}
\begin{document}
\maketitle

\begin{abstract}
In this paper, we present an involution on some kind of colored
$k$-ary trees which provides a combinatorial proof of a
combinatorial sum involving the generalized Catalan numbers
$C_{k,\gamma}(n)=\frac{\gamma}{k n+\gamma}{k n+\gamma\choose n}$.
From the combinatorial sum, we refine the formula for $k$-ary trees
and obtain an implicit formula for the generating function of the
generalized Catalan numbers which obviously implies a Vandermonde
type convolution generalized by Gould. Furthermore, we also obtain a
combinatorial sum involving a vector generalization of the Catalan
numbers by an extension of our involution.

\noindent Mathematics Subject Classification: 05A19, 05C05

\end{abstract}

\section{Introduction}
\label{1} Recently, the author obtained the following identity
involving the Catalan numbers $C_n=\frac{1}{n+1}{2n\choose n}$ by
accident which is similar to the identity $2(a)$ in Riordan's book
\cite[p. 152--153]{6}:
\begin{equation}
\sum_{i=0}^n(-1)^{n-i}{i+1\choose n-i}C_{i}=\delta_{0n}, \mbox{ for
$n\geq 0$},
\end{equation}
where $\delta_{0n}$ is the Kronecker symbol. It is well known
\cite{10} that $C_n$ counts the number of $2$-ary trees (or complete
binary trees) with $n$ internal vertices (vertices with outdegree at
least $1$). Now suppose the number of ordered forests with $\gamma$
$\beta$-ary trees and with total number of $n$ internal vertices is
$C_{\beta,\gamma}(n)$ which is a natural generalization of the
Catalan numbers. It is also well known \cite{4,10} that
$$
C_{\beta,\gamma}(n)=\frac{\gamma}{\beta n+\gamma}{\beta
n+\gamma\choose n}.
$$
Then by generalizing the case $C_n=C_{2,1}(n)$ to
$C_{\beta,\gamma}(n)$, we obtain a generalization identity of $(1)$:
\begin{thm}
\label{1t1}
 For $n\geq
0$, $\alpha,\beta,\gamma\in \mathbb{C}$,
\begin{align}
\sum_{i=0}^n(-1)^{n-i}{(\beta-1)i+\alpha\choose
n-i}\frac{\gamma}{\beta i+\gamma}{\beta i+\gamma\choose
i}=(-1)^n{\alpha-\gamma\choose n}.
\end{align}
\end{thm}
Actually, we can even generalize $(2)$ by introducing following
notations: For any vectors $\vec{a}=(a_1,\ldots,a_t)$ and
$\vec{b}=(b_1,\ldots,b_t)$, we denote $\vec{a}\leq\vec{b} $ if
$a_i\leq b_i$ for all $1\leq i\leq t$;  We also define
$\vec{b}-\vec{a}=(b_1-a_1,\ldots,b_t-a_t)$ and $\vec{a}\cdot
\vec{b}=\sum_{k=1}^ta_kb_k$; As usual, any dimension vector with
constant entries $k$'s will be denoted by $\vec{k}$ for short if
there is no confusion; For $x\in \mathbb{C}$, $\sum_{i=1}^tm_i=x$
and $m_1,\ldots,m_{t-1}\in \mathbb{N}$, we interpret the general
multinomial coefficient as follows:
$$
{x\choose m_1,\ldots,m_{t-1},m_t}={x\choose m_1}{x-m_1\choose
m_2}\cdots{x-m_1-\cdots-m_{t-2}\choose m_{t-1}},
$$
where ${x\choose k}=\frac{x(x-1)\cdots(x-k+1)}{k!}$ for $k\in
\mathbb{N^+}$ and ${x\choose 0}=1$.
\begin{thm}
\label{1t2}
 For $\vec{n}=(n_1,\ldots, n_t)\geq
\vec{0},\vec{p}=(p_1,\ldots,p_t)\geq\vec{1}\in\mathbb{N}^t,t\geq
1,\gamma\geq 0\in\mathbb{N}$ and $\alpha\in\mathbb{C}$, define
$Q(\vec{n};\vec{p};\gamma)$ which can be viewed as a vector
generalization of the Catalan numbers as follows:
\begin{align*}
&Q(\vec{n};\vec{p};\gamma)=\frac{\gamma}{\vec{n}\cdot
\vec{p}+\gamma} {\vec{n}\cdot\vec{p}+\gamma\choose
n_1,\ldots,n_t,\vec{n}\cdot(\vec{p}-\vec{1})+\gamma}.
\end{align*}
Then, there holds
\begin{multline}
\sum_{\vec{i}=(i_1,\ldots,i_t)\geq
\vec{0}}(-1)^{i_1+\cdots+i_t}{(\vec{p}-\vec{1})\cdot
(\vec{n}-\vec{i})+\alpha\choose
i_1,\ldots,i_t,(\vec{p}-\vec{1})\cdot
(\vec{n}-\vec{i})+\alpha-\sum_{j=1}^t
i_j}Q(\vec{n}-\vec{i};\vec{p};\gamma)\\=(-1)^{n_1+\cdots+n_t}{\alpha-\gamma\choose
n_1,\ldots,n_t,\alpha-\gamma-\sum_{j=1}^t n_j}.
\end{multline}
\end{thm}
In fact, $Q(\vec{n};\vec{p};\gamma)$ also counts some kind of
ordered forests which we will say in details later. When looking up
$(2)$ and $(3)$ in the literature, we first found a special case of
$(2)$ for $\alpha=0,\gamma=1$ in a very recent paper of Cohen et al.
\cite[2005]{1}. They claimed that they did not find their result in
the literature and they did not know any obvious combinatorial proof
for their special case, even for $(1)$. But, shortly later, we found
that a special case of $(2)$ for $\alpha=\beta$ also appeared in an
early paper of Gould \cite{3} though he did not state it explicitly.
Hence, the main purpose of this paper is to prove $(2)$ and $(3)$,
especially combinatorially.
\par
Our combinatorial proofs of $(2)$ and $(3)$ are via an involution on
some kind of colored ordered trees. We note that ordered trees and
$k$-ary trees have been extensively studied and refer to \cite{4,10}
and references therein.
In addition, in representing ordered trees we will put the root on
the top and put other vertices on different levels: the root on
level 0, the children of vertices of level $i$ on level $i+1$ (right
below level $i$).
\par
In the rest of this paper, we first present an involution on colored
$k$-ary trees which leads to a combinatorial proof of $(2)$. Next,
we comment briefly on how to prove $(3)$ by a similar argument as
$(2)$. Finally we give a noncombinatorial proof of $(2)$ by a
modified Riordan array theorem derived from theorems in Sprugnoli et
al. \cite{5}. Also we refine the formula for $C_{\beta,\gamma}(n)$
with aids of a Gould classes of inverse relation and obtain an
implicit formula for the generating function of
$C_{\beta,\gamma}(n)$ which obviously implies a Vandermonde type
convolution generalized by Gould \cite{2,3}.

\section{Proofs of $(2),(3)$ via an involution on colored trees}
We rewrite identity $(2)$:
\begin{equation}
\sum_{i=0}^n(-1)^i{(\beta-1)(n-i)+\alpha\choose
i}\frac{\gamma}{\beta (n-i)+\gamma}{\beta (n-i)+\gamma\choose
n-i}=(-1)^n{\alpha-\gamma\choose n}.
\end{equation}
 It is known that the left side of $(4)$ is a polynomial in $\alpha,\beta,\gamma$. Hence it
 suffices to prove $(4)$ for all positive integers
$\alpha,\beta,\gamma$. Therefore, we first give a combinatorial
proof for the case $\gamma=1,\alpha=1,\beta\geq 1$. Then we will
show that with a little more work the general cases will follow.
\par
Firstly, we know that $\frac{1}{\beta(n-i)+1}{\beta(n-i)+1\choose
n-i}$ counts the number of $\beta$-ary trees with $n-i$ internal
vertices and $(\beta-1)(n-i)+1$ is the number of leaves in such a
tree. Hence,
$$
{(\beta-1)(n-i)+1\choose
i}\frac{1}{\beta(n-i)+1}{\beta(n-i)+1\choose n-i}
$$
counts the number of $\beta$-ary trees with $n-i$ internal vertices
and $i$ colored leaves. (Note: When a $\beta$-ary tree has only one
vertex, this vertex will be treated as a leaf.)
\par
However, we can classify all those $\beta$-ary trees with $n-i$
internal vertices and $i$ colored leaves into two classes for
$n>0,i\geq 0$:
\begin{itemize}
 \item The first
class consists of trees in which there is at least one colored leaf
in the lowest two levels (the two levels furthest to the root) and
to the left of which all vertices in the same level have no
children. The leftmost and lowest such a colored leaf will be called
the \emph{candidate leaf}.
\item
The second class consists of the rest of trees, i.e., the trees in
which there is no colored leaf in the lowest level (the level
furthest to the root) and there is an internal vertex in the second
lowest level to whose left in the same level there is no colored
leaf. The leftmost such an internal vertex will be called the
\emph{incumbent internal vertex}.
\end{itemize}

\par
Now, we define an operation on these colored $\beta$-ary trees which
will give an involution:
\begin{itemize}
\item
For each $\beta$-ary tree with $n-i$ internal vertices and $i$
colored leaves of the first class, we attach $\beta$ leaves to its
candidate leaf and erase the color of the candidate leaf (now an
internal vertex). Thus we obtain a $\beta$-ary tree with $n-i+1$
internal vertices and $i-1$ colored leaves.
\item
For each $\beta$-ary tree of the second class, we delete all
children of its incumbent internal vertex and color the incumbent
internal vertex (now a leaf). Thus, we obtain a $\beta$-ary tree
with $n-i-1$ internal vertices and $i+1$ colored leaves.
\end{itemize}
Figure 1 shows an example of this involution. From above involution,
when $n>0$, we obtain a perfect matching over all $\beta$-ary trees
with $n-i$ internal vertices and $i$ colored leaves for all $i\geq
0$. Moreover, if each colored ordered tree $T$ has weight
$(-1)^{\#\mbox{\scriptsize of colored leaves in $T$}}$, then the two
trees matched to each other are weighted by $+1$ and $-1$
respectively. Hence, the total weight over these colored $\beta$-ary
trees is 0. The case for $n=0$ is obvious. This completes the proof
of the case $\alpha=1,\gamma=1$ of $(2)$
$$
\sum_{i=0}^{n}(-1)^i{(\beta-1)(n-i)+1\choose
i}\frac{1}{\beta(n-i)+1}{\beta(n-i)+1\choose n-i}=\delta_{0n}.\mbox{
\qed}
 $$
\par
\begin{center}
\setlength{\unitlength}{1mm}
\begin{picture}(100,40)
\put(20,40){\circle{1.5}}\put(20,40){\line(-1,-1){10}}
\put(20,40){\line(0,-1){10}}\put(20,40){\line(1,-1){10}}
\multiput(10,30)(20,0){2}{\circle{1.5}}\put(20,30){\circle*{1.5}}
\put(10,30){\line(-1,-1){10}}
\put(10,30){\line(-2,-5){4}}\put(10,30){\line(2,-5){4}}
\multiput(0,20)(6,0){2}{\circle{1.5}}\put(14,20){\circle{1.5}}
\put(30,30){\line(-1,-1){10}}
\put(30,30){\line(-2,-5){4}}\put(30,30){\line(3,-5){6}}
\put(20,20){\circle{1.5}}\put(26,20){\circle*{1.5}}\put(36,20){\circle{1.5}}
\put(24.5,18.5){$\square$}
\put(36,20){\line(-2,-5){4}}\put(36,20){\line(2,-5){4}}\put(36,20){\line(1,-1){10}}
\put(32,10){\circle{1.5}}\put(40,10){\circle{1.5}}\put(46,10){\circle{1.5}}
\put(74,40){\circle{1.5}}\put(74,40){\line(-1,-1){10}}
\put(74,40){\line(0,-1){10}}\put(74,40){\line(1,-1){10}}
\multiput(64,30)(20,0){2}{\circle{1.5}}\put(74,30){\circle*{1.5}}
\put(64,30){\line(-1,-1){10}}
\put(64,30){\line(-2,-5){4}}\put(64,30){\line(2,-5){4}}
\multiput(54,20)(6,0){2}{\circle{1.5}}\put(68,20){\circle{1.5}}
\put(84,30){\line(-1,-1){10}}
\put(84,30){\line(-2,-5){4}}\put(84,30){\line(3,-5){6}}
\put(74,20){\circle{1.5}}\put(80,20){\circle{1.5}}\put(90,20){\circle{1.5}}
{\thicklines\put(80,20){\circle{2.5}}}
\put(90,20){\line(-2,-5){4}}\put(90,20){\line(2,-5){4}}\put(90,20){\line(1,-1){10}}
\put(86,10){\circle{1.5}}\put(94,10){\circle{1.5}}\put(100,10){\circle{1.5}}
\put(80,20){\line(-1,-1){10}}\put(80,20){\line(-2,-5){4}}\put(80,20){\line(2,-5){4}}
\put(70,10){\circle{1.5}}\put(76,10){\circle{1.5}}\put(84,10){\circle{1.5}}
\put(43,30){$\Longleftrightarrow$}
\put(1,5){\circle*{1.5} }\put(-10,-2){\text{a colored
leaf}}\put(30,5){\circle*{1.5}}\put(28.5,3.5){$\square$}\put(20,-2){\text{the
candidate leaf
}}\put(66,5){\circle{1.5}}{\thicklines\put(66,5){\circle{2.5}}}\put(60,-2){\text{the
incumbent internal vertex}}
\end{picture}
\end{center}
\begin{center}
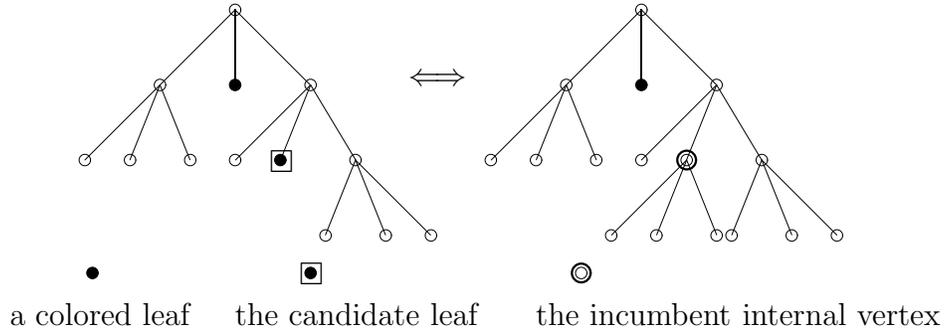
Figure 1: An example of above involution for $n=6,\alpha=1,\beta=3$,
the tree on the left side is of the first class while the other
belongs to the second class.
\end{center}

For simplicity of the remaining presentment, we describe the core
idea of the above involution in the following obvious lemma.
\begin{lem}\label{3l1}
Suppose $S$ is a finite set of structures each of which has two
finite integer indices (or called parameters, etc.), say $ x$ and
$y$, and has weight $(-1)^x$. All these structures can be
partitioned into two classes: the first and the second. Suppose
every structure of the first class in $S$ corresponds uniquely to a
structure of the second class in $S$ with index $x$ increasing by
$1$ and index $y$ decreasing by $1$ and vice versa. Then the total
weight of structures in $S$ is 0.
\end{lem}

 Now, when
$\alpha>1$, we just view
$$
{(\beta-1)(n-i)+\alpha\choose
i}\frac{1}{\beta(n-i)+1}{\beta(n-i)+1\choose n-i}
$$
as the number of $(\alpha-1)$-planted $\beta$-ary trees (analog to
usual planted tree except that there are $\alpha-1$ planted roots,
see Figure 2 for an example) with $n-i$ internal vertices and $j$
colored leaves and $i-j$ colored planted roots for all $0\leq
i-j\leq \alpha-1$. (Note: The planted roots are neither internal
vertices nor leaves.)
\par
\begin{center}
\setlength{\unitlength}{1mm}
\begin{picture}(60,40)
\put(20,20){\line(-3,2){15}}\put(20,20){\line(-6,5){12}}\put(20,20){\line(-1,1){10}}
\put(20,20){\line(-3,5){6}}\put(20,20){\line(3,5){6}}\put(20,20){\line(1,1){10}}
\multiput(5,30)(3,0){2}{\circle{1.5}}\put(10,30){\circle*{1.5}}\put(14,30){\circle{1.5}}
\put(26,30){\circle*{1.5}}\put(30,30){\circle{1.5}}
\put(18,29){\text{$\cdots$}}\put(5,32){\text{$\overbrace{\hskip
2.5cm}$}}\put(3,36){\text{$\alpha-1$ planted roots}}

\put(20,20){\circle{1.5}}\put(20,20){\line(-1,-1){10}}
\put(20,20){\line(0,-1){10}}\put(20,20){\line(1,-1){10}}
\multiput(10,10)(20,0){2}{\circle{1.5}}\put(20,10){\circle*{1.5}}
\put(10,10){\line(-1,-1){10}}
\put(10,10){\line(-2,-5){4}}\put(10,10){\line(2,-5){4}}
\multiput(0,0)(6,0){2}{\circle{1.5}}\put(14,0){\circle{1.5}}
\put(30,10){\line(-1,-1){10}}
\put(30,10){\line(-2,-5){4}}\put(30,10){\line(3,-5){6}}
\put(20,0){\circle{1.5}}\put(26,0){\circle*{1.5}}\put(36,0){\circle{1.5}}
\put(24.5,-1.5){$\square$}
\end{picture}
\end{center}
\begin{center}
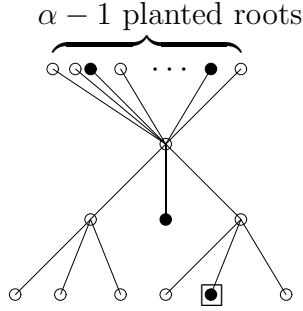
Figure 2: An example of $(\alpha-1)$-planted $3$-ary tree \\in which
some planted roots are also colored.
\end{center}
Following Lemma \ref{3l1}, the two indices here are respectively the
number of internal vertices ($y$) and the sum of colored leaves
together with the number of colored planted roots ($x$), and the two
classes are defined to be the same as in the above proof for
$\alpha=1,\gamma=1$.
\begin{itemize}
\item[(i)]
When $n>\alpha-1>0$, for all possible $0\leq i< n$, any planted
$\beta$-ary tree with $n-i$ internal vertices and $j$ colored leaves
and $i-j$ colored planted roots is either of the first class or of
the second class. So their total weight is 0 from our involution.
\item[(ii)]
When $n\leq \alpha-1$, except for those structures which belong to
the defined two classes, there still remain structures for $i=n$ and
with only the planted roots colored. See Figure 3 for an example of
the case $i=n$. Since the total weight of structures belonging to
the defined two class is 0, we only need to calculate the total
weight of those exceptional structures. This is equivalent to
choosing $n$ out of all $\alpha-1$ planted roots. Hence, the total
weight is
$
(-1)^n{\alpha-1\choose n}.
$
\end{itemize}
Thus the proof for all $\alpha\geq 1,\beta\geq 1$ follows.\qed
\begin{center}
\setlength{\unitlength}{1mm}
\begin{picture}(60,20)
\put(20,0){\line(-3,2){15}}\put(20,0){\line(-6,5){12}}\put(20,0){\line(-1,1){10}}
\put(20,0){\line(-3,5){6}}\put(20,0){\line(3,5){6}}\put(20,0){\line(1,1){10}}
\multiput(5,10)(3,0){2}{\circle{1.5}}\put(10,10){\circle*{1.5}}\put(14,10){\circle{1.5}}
\put(26,10){\circle*{1.5}}\put(30,10){\circle{1.5}}
\put(18,9){\text{$\cdots$}}\put(5,12){\text{$\overbrace{\hskip
2.5cm}$}}\put(3,16){\text{$\alpha-1$ planted roots}}

\put(20,0){\circle{1.5}}
\end{picture}
\end{center}
\begin{center}
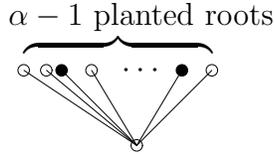
Figure 3: An example of $(\alpha-1)$-planted $\beta$-ary tree with
$0$ \\internal vertices and $n$ colored plated roots, i.e., $i=n$.
\end{center}
%
\par
For $\gamma\geq 1$, it is well known that
$$
\frac{\gamma}{\beta n+\gamma}{\beta n+\gamma\choose n}
$$
counts the number of ordered forests of $\beta$-ary trees with
totally $n$ internal vertices and $\gamma$ components and in which
there are just $(\beta-1)n+\gamma$ leaves. Thus, we can view
$$
{(\beta-1)(n-i)+\alpha\choose
i}\frac{\gamma}{\beta(n-i)+\gamma}{\beta(n-i)+\gamma\choose n-i}
$$
as the number of $(\alpha-\gamma)$-planted $\beta$-ary forests with
$\gamma$ components where the planted roots are planted on the first
ordered tree in the forest. We put all roots (except planted roots)
of the forest on the same level and define the two classes of those
forests in the same way as trees.
 So, with same
argument, the general case $\gamma\geq 1$ follows. See an
illustration in Figure 4. Thus, we have completed a combinatorial
proof of $(2)$ for all positive integers $\alpha,\beta,\gamma$,
which implies that $(2)$ holds for all
$\alpha,\beta,\gamma\in\mathbb{C}$.\qed
\par
\vskip 15pt
\begin{center}
\setlength{\unitlength}{1mm}
\begin{picture}(120,40)
\put(20,30){\circle{1.5}}\put(20,30){\line(-1,4){2.5}}\put(20,30){\line(0,1){10}}
\put(20,30){\line(1,4){2.5}}\multiput(17.5,40)(5,0){2}{\circle{1.5}}\put(20,40){\circle*{1.5}}
\put(20,30){\line(-1,-1){10}}\put(20,30){\line(-2,-5){4}}\put(20,30){\line(2,-5){4}}
\put(10,20){\circle{1.5}}\put(16,20){\circle{1.5}}\put(24,20){\circle*{1.5}}
\put(10,20){\line(-1,-1){10}}\put(10,20){\line(-3,-5){6}}\put(10,20){\line(0,-1){10}}
\put(0,10){\circle{1.5}}\put(4,10){\circle{1.5}}\put(10,10){\circle{1.5}}
\put(16,20){\line(-2,-5){4}}\put(16,20){\line(0,-1){10}}\put(16,20){\line(3,-5){6}}
\put(12,10){\circle{1.5}}\put(16,10){\circle{1.5}}\put(22,10){\circle{1.5}}
\put(30,30){\circle{1.5}}\put(30,30){\line(-2,-5){4}}\put(30,30){\line(1,-5){2}}
\put(30,30){\line(4,-5){8}}
\put(26,20){\circle{1.5}}\put(32,20){\circle{1.5}}\put(38,20){\circle*{1.5}}
\put(32,20){\line(-2,-5){4}}\put(32,20){\line(1,-5){2}}\put(32,20){\line(4,-5){8}}
\put(28,10){\circle*{1.5}}\put(26.5,8.5){$\square$}\put(34,10){\circle{1.5}}\put(40,10){\circle{1.5}}
\put(34,10){\line(-2,-5){4}}\put(34,10){\line(1,-5){2}}\put(34,10){\line(4,-5){8}}
\put(30,0){\circle{1.5}}\put(36,0){\circle{1.5}}\put(42,0){\circle{1.5}}


\put(70,30){\circle{1.5}}\put(70,30){\line(-1,4){2.5}}\put(70,30){\line(0,1){10}}
\put(70,30){\line(1,4){2.5}}\multiput(67.5,40)(5,0){2}{\circle{1.5}}\put(70,40){\circle*{1.5}}
\put(70,30){\line(-1,-1){10}}\put(70,30){\line(-2,-5){4}}\put(70,30){\line(2,-5){4}}
\put(60,20){\circle{1.5}}\put(66,20){\circle{1.5}}\put(74,20){\circle*{1.5}}
\put(60,20){\line(-1,-1){10}}\put(60,20){\line(-3,-5){6}}\put(60,20){\line(0,-1){10}}
\put(50,10){\circle{1.5}}\put(54,10){\circle{1.5}}\put(60,10){\circle{1.5}}
\put(66,20){\line(-2,-5){4}}\put(66,20){\line(0,-1){10}}\put(66,20){\line(3,-5){6}}
\put(62,10){\circle{1.5}}\put(66,10){\circle{1.5}}\put(72,10){\circle{1.5}}
\put(80,30){\circle{1.5}}\put(80,30){\line(-2,-5){4}}\put(80,30){\line(1,-5){2}}
\put(80,30){\line(4,-5){8}}
\put(76,20){\circle{1.5}}\put(82,20){\circle{1.5}}\put(88,20){\circle*{1.5}}
\put(82,20){\line(-2,-5){4}}\put(82,20){\line(1,-5){2}}\put(82,20){\line(4,-5){8}}
\put(78,10){\circle{1.5}}\put(78,10){\circle{2.5}}\put(84,10){\circle{1.5}}\put(90,10){\circle{1.5}}
\put(84,10){\line(-1,-5){2}}\put(84,10){\line(1,-5){2}}\put(84,10){\line(4,-5){8}}
\put(82,0){\circle{1.5}}\put(86,0){\circle{1.5}}\put(92,0){\circle{1.5}}
\put(78,10){\line(-1,-1){10}}\put(78,10){\line(-2,-5){4}}\put(78,10){\line(1,-5){2}}
\put(45,25){$\Longleftrightarrow$}
\put(68,0){\circle{1.5}}\put(74,0){\circle{1.5}}\put(80,0){\circle{1.5}}

\end{picture}
\end{center}
\begin{center}
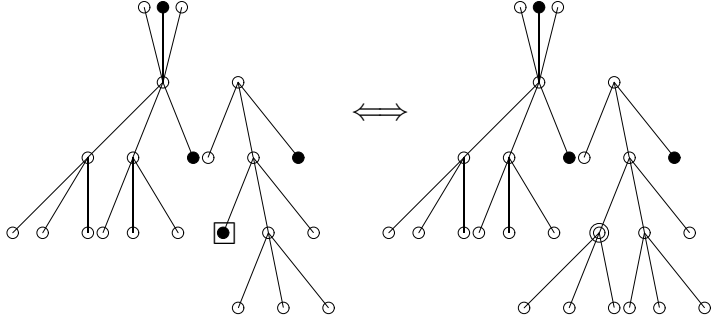
Figure 4: An example of the involution for
$\alpha=4,\beta=3,\gamma=2,n=10$.
\end{center}

To prove $(3)$, we first note the following facts:
\begin{fact}\label{3f2}(Stanley \cite[p.34]{10})
The number of ordered forests with $n_i$ internal vertices of
outdegree $p_i$ and $\gamma$ component trees is
$$
Q(\vec{n};\vec{p};\gamma)=\frac{\gamma}{\vec{n}\cdot\vec{p}+\gamma}{\vec{n}\cdot\vec{p}+\gamma\choose
n_1,\ldots,n_t,\vec{n}\cdot\vec{p}+\gamma-\sum_{j=1}^tn_j}.
$$
In addition, the number of leaves in such a forest is
$\vec{n}\cdot(\vec{p}-\vec{1})+\gamma$.
\end{fact}
Thus, the summand in the left hand side of $(3)$ counts the total
weight of $(\alpha-\gamma)$-planted ordered forests with leaves (and
planted roots) colored with colors in $\{1,\ldots,t\}$, and in which
there are $n_j-i_j$ internal vertices of outdegree $p_j$ and $i_j$
leaves plus planted roots have color ${j}$ for $1\leq j\leq t$. We
also assume that all planted roots are planted on the first ordered
tree when $r>1$ and the two classes are still the same. But we
modify the involution a little bit:
\begin{itemize}
\item[(i)]
 If the candidate leaf has color $j$,
then we attach $p_j$ leaves to it and erase its color.
\item[(ii)]
Correspondingly, if the incumbent internal vertex has $p_j$ leaves,
we delete all its leaves and color it with color $j$.
\end{itemize}
We omit figures for illustration of this case since it is so clear
and simple. After applying this modified involution and by the same
reasoning as ordered $\beta$-ary forests, all cancel out except for
the case $i_j=n_j$ with only planted roots colored. It reduces to
choosing out $n_i$ planted roots from $\alpha-\gamma$ and coloring
them with color $i$ for all $1\leq i\leq t$, for which there are
$$
{\alpha-\gamma\choose n_1,\ldots,n_t,\alpha-\gamma-\sum_{i=1}^tn_i}
$$
different ways. Thus we obtain $(3)$. However, because the left side
of $(3)$ is a polynomial in $\alpha$ but not necessarily a
polynomial in $\vec{p}$ and $\gamma$, we can not generalize it to
$\vec{p}\in\mathbb{C}^t$ and $\gamma\in \mathbb{C}$ in the manner of
$(2)$.
\section{Refinement of $C_{\beta,\gamma}(n)$ and the Gould-Vandermonde's convolution}
\label{2}

 In this section, we will
prove $(2)$ and refine the formula for $k$-ary trees and the
Gould-Vandermonde's convolution through Riordan array theory. The
Riordan array theory is an effective tool to prove and invert
combinatorial identities \cite{5,6,7,8,9}. For clarity, here we
present the related results which will be used later.
\begin{defn}
A lower triangular infinite matrix, $T$, is a Riordan array, if the
generating function of the $k$th column is
$$
g(x)f(x)^k,
$$
for $k=0,1,2,3,\ldots$, where
\begin{align*}
g(x)&=g_0+g_1x+g_2x^2+g_3x^3+\cdots,\quad g_0\neq 0,\\
f(x)&=\qquad f_1x+f_2x^2+f_3x^3+\cdots, \quad f_1\neq 0.
\end{align*}
In addition, we will denote this Riordan array as $[g(x),f(x)]$ and
denote its $(n,k)$ entry as $[g(x),f(x)]_{n,k}$.
\end{defn}
\begin{thm}[Riordan array theorem]\label{2t1}
Let $A(x)=a_0+a_1x+a_2x^2+\cdots$, and
$L(x)=l_0+l_1x+l_2x^2+l_3x^3+\cdots$, and $[g(x),f(x)]$ is a Riordan
array. Then
$$
\sum_{k=0}^n[g(x),f(x)]_{n,k}a_k=l_n, \mbox{\quad for
$n=0,1,2,\ldots$},
$$
if and only if
\begin{equation}
g(x)A(f(x))=L(x).
\end{equation}
\end{thm}
Sprugnoli et al. \cite{5} studied how to invert combinatorial sums
through Riordan array theory and Lagrange inversion formula. We
present the following theorem easily derived from Theorem 2.1--2.4
in the paper \cite{5}.
\begin{thm}[modified Riordan array theorem]\label{2t2}
Let $A(x)=a_0+a_1x+a_2x^2+\cdots$, and
$L(x)=l_0+l_1x+l_2x^2+l_3x^3+\cdots$, and $[g(x),f(x)]$ is a Riordan
array. Then
$$
\sum_{k=0}^n[g(x),f(x)]_{n,k}a_k=l_n, \mbox{\quad for
$n=0,1,2,\ldots$},
$$
if and only if
\begin{equation}
n[x^n]A(x)=[x^{n-1}]\Big(\frac{x}{f(x)}\Big)^n\Big(\frac{L(x)}{g(x)}\Big)',
\quad\mbox{for $n>0$},
\end{equation}
and $a_0=\frac{L(0)}{g(0)}$.

\end{thm}

Now, we come to a proof of $(2)$ via the modified Riordan array
theorem:
Firstly, for $(2)$ we have $A(x)=\sum_{n\geq 0}\frac{\gamma}{\beta
n+\gamma}{\beta n+\gamma\choose n}x^n,
f(x)=x(1-x)^{\beta-1},g(x)=(1-x)^{\alpha},L(x)=\sum_{n\geq
0}(-1)^n{\alpha-\gamma\choose n}x^n=(1-x)^{\alpha-\gamma}$. Thus,
for $n\geq 1$,
\begin{align*}
[x^{n-1}]\Big(\frac{x}{x(1-x)^{\beta-1}}\Big)^n\Big(\frac{(1-x)^{\alpha-\gamma}}{(1-x)^\alpha}\Big)'&=[x^{n-1}]\gamma(1-x)^{(1-\beta)n-\gamma-1}\\
&={\gamma}{\beta n+\gamma-1\choose n-1}\\
&=n[x^n]\sum_{n\geq 0}\frac{\gamma}{\beta n+\gamma}{\beta
n+\gamma\choose n}x^n.
\end{align*}
Additionally, $\frac{\gamma}{\beta 0+\gamma}{\beta 0+\gamma\choose
0}=\frac{(1-0)^{\alpha-\gamma}}{(1-0)^\alpha}$. By Theorem
\ref{2t2}, the proof follows.\qed
\par
Let $C_{\beta,\gamma}(x)=\sum_{n\geq 0}C_{\beta,\gamma}(n)x^n$. By
Riordan array theorem $(5)$ and $(2)$, we also obtain an implicit
formula for $C_{\beta,\gamma}(x)$ which appeared in Gould \cite{3}.
\begin{cor}\label{2c1}
\begin{equation}
C_{\beta,\gamma}(x(1-x)^{\beta-1})=(1-x)^{-\gamma}.
\end{equation}
\end{cor}
From $(7)$, we immediately obtain
$C_{\beta,\alpha_1+\alpha_2}(x)=C_{\beta,\alpha_1}(x)C_{\beta,\alpha_2}(x)$
which leads to the following famous convolution of Vandermonde type
generalized by Gould \cite{2,3}. (We refer to it as
Gould-Vandermonde's convolution in the title of this paper.)
\begin{cor}\label{2c2} For any $\alpha,\beta,\gamma\in \mathbb{C}$
and $n\geq 0$, there holds
\begin{multline}
\frac{\alpha_1+\alpha_2}{\beta n+\alpha_1+\alpha_2}{\beta
n+\alpha_1+\alpha_2\choose n}\\
=\sum_{i=0}^n\frac{\alpha_1}{\beta i+\alpha_1}{\beta
i+\alpha_1\choose i}\frac{\alpha_2}{\beta (n-i)+\alpha_2}{\beta
(n-i)+\alpha_2\choose n-i}.
\end{multline}
\end{cor}

Note that even if we did not know the exact formula of
$C_{\beta,\gamma}(n)$, we still have by the involution argument in
section 2 that
$$
\sum_{i=0}^n(-1)^{n-i}{(\beta-1)i+\alpha\choose
n-i}C_{\beta,\gamma}(i)=(-1)^n{\alpha-\gamma\choose n}.
$$
Hence, from the Gould classes of inverse relation \cite[p.52]{6}
\begin{align}
b_n=\sum_{k=0}^n{m+ak\choose n-k}z^{n-k}a_k,
a_n=\sum_{k=0}^{n}\frac{-ak-m}{-an-m}{-an-m\choose n-k}b_kz^{n-k},
\end{align}
 we obtain
\begin{cor} \label{2c3}
For $n>0$,
\begin{align}
C_{\beta,\gamma}(n)=\sum_{k=0}^{n}(-1)^n\frac{(1-\beta)k-\alpha}{(1-\beta)n-\alpha}{(1-\beta)n-\alpha\choose
n-k}{\alpha-\gamma\choose k}.
\end{align}

\end{cor}
By setting $\alpha=0$ in $(10)$ and applying Vandermonde's
convolution, it is easy to obtain a closed formula:
\begin{align*}
C_{\beta,\gamma}(n)&=\sum_{k=0}^n(-1)^n\frac{n-(n-k)}{n}{(1-\beta) n\choose n-k}{-\gamma\choose k}\\
 &=\sum_{k=0}^n(-1)^n{(1-\beta) n\choose n-k}{-\gamma\choose k}-\sum_{k=0}^n(-1)^n\frac{n-k}{n}{(1-\beta) n\choose n-k}{-\gamma\choose k}\\
&=(-1)^n{(1-\beta) n-\gamma\choose n}+(-1)^n(\beta-1){(1-\beta)
n-1-\gamma\choose n-1}\\
 &=\frac{\gamma}{\beta n+\gamma}{\beta
n+\gamma\choose n}.
\end{align*}

\par
We finally remark that our approach to obtain the formula for
$k$-ary trees may be applied to other combinatorial structures,
especially recursive ones.

\vskip 20pt

 {\noindent\bf Acknowledgement.} The author is very grateful to the referees for their
 valuable suggestions and comments. 

\end{document}